\documentclass[12pt]{amsart}
\usepackage{amsfonts}
\usepackage{amssymb}
\theoremstyle{definition}

\theoremstyle{remark}

\newcommand{\ds}{\displaystyle}

\begin{document}

\centerline{\large\bf CONFORMAL \ FLAT \ $AK_2$-MANIFOLDS
\footnote{\it PLISKA Studia mathematica bulgarica. Vol. 9, 1987, p. 12-16.}}

\vspace{0.2in}
\centerline{\large OGNIAN \ T.  \ KASSABOV}

\vspace{0.3in}
{\sl In this note we examine $AK_2$-manifolds of dimension $2m \ge 6$.}

\vspace{0.1in}
{\bf 1. Introduction.} Let $M$ be a $2m$ - dimensional almost Hermitian manifolds with
metric $g$ and almost complex structure $J$. The Riemannian connection and the curvature
tensor are denoted by $\nabla$ and $R$, respectively. The manifold is said to be a
K\"ahler, or nearly K\"ahler, or almost K\"ahler manifold, if
$$
	\begin{array}{l}
		\nabla J=0\quad {\rm or} \quad (\nabla_XJ)X=0, \quad {\rm or} \\
		 \\
		g((\nabla_XJ))Y,Z)+g((\nabla_YJ))Z,X)+g((\nabla_ZJ))X,Y)=0 \ ,
	\end{array}   \leqno (1.1)
$$
respectively. The corresponding classes of manifolds are denoted by $K$, $NK$, $AK$,
respectively. It is well known, that for these classes
$$
	(\nabla_XJ)Y + (\nabla_{JX}J)JY=0     \leqno (1.2)
$$
holds [2].

For a given class $L$ of almost Hermitian manifolds its subclass $L_i$ is defined by the
identity $(i)$, where

1) $R(X,Y,Z,U)=R(JX,JY,Z,U)$,

2) $R(X,Y,Z,U)=R(JX,JY,Z,U)+R(JX,Y,JZ,U)+R(JX,Y,Z,JU)$,

3) $R(X,Y,Z,U)=R(JX,JY,JZ,JU)$.

\noindent
Then we have $L_1 \subset L_2 \subset L_3$ and $NK =NK_2$, $K=NK_1=AK_1$, $K=NK\cap AK$ [2].

The Weil conformal curvature tensor $C$ for $M$ is defined by
$$
	\begin{array}{r}
		\ds C(X,Y,Z,U)=R(X,Y,Z,U)- \frac{1}{2m-2} \{ g(X,U)S(Y,Z) \\
		 \\
		-g(X,Z)S(Y,U)+g(Y,Z)S(X,U)-g(Y,U)S(X,Z) \} \\ 
		\\
		\ds +\frac{\tau}{(2m-1)(2m-2)} \{ g(X,U)g(Y,Z) - g(X,Z)g(Y,U) \} \ ,
	\end{array}   
$$
where $S$ and $\tau$ are the Ricci tensor and the scalar curvature, respectively. It is
well known, that (if $m\ge 2$) $M$ is conformal flat if and only if  $C=0$.

Conformal K\"ahler and nearly K\"ahler manifolds are classified in [4] and [5]. Here,
we shall prove the following theorem:

T h e o r e m. {\it Let $M \in AK_2$ be a $2m$-dimensional conformal flat manifold, 
$m>2$. Then it is one of the following:

a) a flat K\"ahler manifold;

b) a 6-dimensional almost K\"ahler manifold of constant negative sectional curvature;

c) locally $M_1 \times M_2$, where $M_1$ (resp. $M_2$) is a 4-dimensional almost
K\"ahler manifold of constant sectional curvature $-c$ (resp. a 2-dimensional
K\"ahler manifold of constant sectional curvature $c$), $c>0$;

d) locally $M_3 \times M_2$, where $M_3$ is a  6-dimensional almost K\"ahler manifold 
of constant negative sectional curvature $-c$.}

R e m a r k \ 1. We don't know whether there exists an almost K\"ahler manifold of 
constant negative sectional curvature of dimension 4 or 6.

R e m a r k \ 2. If a conformal flat almost Hermitian manifold $M$ satisfies the
identity 3), then $S(X,Y)=S(JX,JY)$ and $M$ satisfies also the identity 2).

\vspace{0.1in}
{\bf 2. Preliminaries.} Let $Q$ be a tensor of type (1.1). According to the Ricci
identity,
$$
	(\nabla_X(\nabla_YQ))Z-(\nabla_Y(\nabla_XQ))Z = 
	R(X,Y)QZ - QR(X,Y)Z \ .     \leqno (2.1)
$$ 
From the second Bianchi identity it follows
$$
	\sum_{i=1}^{2m} (\nabla_{E_i}R)(X,Y,Z,E_i)=
	(\nabla_XS)(Y,Z) - (\nabla_YS)(X,Z) \ ,  \leqno (2.2)
$$
$$
	\sum_{i=1}^{2m} (\nabla_{E_i}S)(X,E_i)=	\frac 12 X(\tau) \ ,  \leqno (2.3)
$$
where $\{ E_i, \ i=1,...,2m \}$ is a local orthonormal frame field. We shall assume that
$ E_{m+i}=JE_i$, $i=1,...,m$.

Let the tensor $S'$ be defined by 
$$
	S'(X,Y) = \sum_{i=1}^{2m} R(X,E_i,JE_i,JY)\ .
$$

For an $AK_2$-manifold the following identities [1,2] hold:
$$
	2(\nabla_X(S-S'))(Y,Z) = (S-S')((\nabla_XJ)Y,JZ)+(S-S')(JY,(\nabla_XJ)Z) \ , \leqno (2.4)
$$
$$
	\sum_{i=1}^{2m} (\nabla_{E_i}(\nabla_{E_i}J))Y = \sum_{i=1}^{2m} J(\nabla_{E_i}J)(\nabla_{E_i}J)Y\, , \ \leqno(2.5)
$$
$$
	R(X,Y,Z,U)-R(X,Y,JZ,JU)=\frac 12 g(K(X,Y),K(Z,U))\ , \leqno (2.6)
$$
where $K(X,Y) = (\nabla_XJ)Y-(\nabla_YJ)X$.

A 2-dimensional almost Hermitian manifold is a K\"ahler manifold. It follows easily from (2.6), 
that if $M$ is an almost K\"ahler manifold of constant curvature $c$ and if $ \dim M \ge 4$, 
then $c \le 0$ and $c=0$ if and only if $M$ is a K\"ahler manifold. On the other hand, an
almost K\"ahler manifold of constant sectional curvature and dimension $2m \ge 8$ is automatically a K\"ahler manifold [3].

\vspace{0.1in}
{\bf 3. Proof of the theorem.} From $C=0$, (2.2) and (2.3) it follows
$$
	(\nabla_XS)(Y,Z)-(\nabla_YS)(X,Z)=\ds \frac {1}{2(2m-1)}
	\{ g(Y,Z)X(\tau) - g(X,Z)Y(\tau) \} \ .  \leqno (3.1) 
$$

Since $C=0$ we have
$$
	S'= \frac{1}{m-1}S-\frac{\tau}{2(m-1)(2m-1)}g \ .
$$
Hence, using (2.4), we find
$$
	2(\nabla_XS)(Y,Z) = S((\nabla_XJ)Y,JZ)+S(JY,(\nabla_XJ)Z) 
	-\frac{X(\tau)}{(m-1)(2m-1)}g(Y,Z)   \ . \leqno (3.2)
$$
Let $X \perp Y,\ JY$. According to $(3.2)$ and (1.2),
$$
	(\nabla_XS)(Y,Y)+(\nabla_XS)(JY,JY) - (\nabla_YS)(X,Y)-(\nabla_{JY}S)(X,JY)
	=-\frac{X(\tau)g(Y,Y)}{(m-1)(2m-1)} \ .
$$ 
The last equality and (3.1) give $X(\tau)=0$. From $X(\tau)=0$, (3.1) and (3.2)
we obtain
$$
	(\nabla_XS)(Y,Z)=(\nabla_YS)(X,Z) \ ,  \leqno (3.3)
$$
$$
	2(\nabla_XS)(Y,Z) = S((\nabla_XJ)Y,JZ)+S(JY,(\nabla_XJ)Z) 
	   \ . \leqno (3.4)
$$

Now let $p \in M$ and $ \{ e_i;\, i=1,...,2m \}$ be an orthonormal basis of
$T_p(M)$, such that $ e_{i+m}=Je_i$ and $Se_i=\lambda e_i$ for $i=1,...,m$.
Let $ \{ E_i; \, i=1,...,2m \} $ be a local orthonormal frame field, such
that $ E_i|p = e_i$ for $i=1,...,2m$. We have
$$
	\sum_{i=1}^{2m} (\nabla_{e_i}(\nabla_{e_i}S))(e_j,e_j)=	
$$
$$
	=\sum_{i=1}^{2m} \{ (\nabla_{E_i}(\nabla_{E_i}S))(E_j,E_j)
	-(\nabla_{\nabla_{E_i}E_i}S)(E_j,E_j)
	-2(\nabla_{E_i}S)(\nabla_{E_i}E_j,E_j) \}_p \ \ \ {\rm using \ (3.4)}
$$
$$
	=\sum_{i=1}^{2m} \{ (\nabla_{E_i}S)((\nabla_{E_i}J)E_j,JE_j)+
	S((\nabla_{E_i}(\nabla_{E_i}J)E_j,JE_j) 
$$
$$
	+S((\nabla_{E_i}J)E_j,(\nabla_{E_i}J)E_j) \}_p \ \ \ {\rm using \ (2.5) \ and \ (3.4) }	
$$
$$
	=-\sum_{i=1}^{2m} (\nabla_{e_i}S)((\nabla_{e_i}J)e_j,Je_j) \ \ \ {\rm using \ (3.3) }
$$
$$
	=\sum_{i=1}^{2m} (\nabla_{(\nabla_{e_i}J)e_j} S)(e_i,Je_j), \ \ \ {\rm using (3.4) }
$$
$$
	=\frac 12 \sum_{i=1}^{2m} (\lambda_j - \lambda_i) g((\nabla_{(\nabla_{e_i}J)e_j}J)e_i,e_j)
$$
and using (1.1), we obtain
$$
	\begin{array}{c}
	\ds\sum_{i=1}^{2m} (\nabla_{e_i}(\nabla_{e_i}S))(e_j,e_j) \\
	\ds =\frac 12 \sum_{i=1}^{2m} (\lambda_j-\lambda_i)
	\{ g((\nabla_{e_i}J)e_j,(\nabla_{e_j}J)e_i) - g(((\nabla_{e_i}J)e_j,(\nabla_{e_i}J)e_j) \}  . 
	\end{array}   \leqno(3.5)
$$

Because of $X(\tau)=0$ and (3.3) we have
$$
	\sum_{i=1}^{2m} (\nabla_{E_j}(\nabla_{E_i}S))(E_i,E_j) =  0 \ .
$$
Using (3.3), we obtain also
$$
	\sum_{i=1}^{2m} (\nabla_{E_i}(\nabla_{E_i}S))(E_j,E_j)=
	\sum_{i=1}^{2m} (\nabla_{E_i}(\nabla_{E_j}S))(E_i,E_j) \ .
$$
From the last two equalities and (2.1) it follows
$$
	\sum_{i=1}^{2m} (\nabla_{e_i}(\nabla_{e_i}S))(e_j,e_j)=
	\sum_{i=1}^{2m} (\lambda_j-\lambda_i)R(e_i,e_j,e_j,e_i) \ .  \leqno (3.6)
$$

Now we compute
$$
	  (\nabla_{e_i}(\nabla_{e_j}S))(e_i,e_j)-
	  (\nabla_{e_j}(\nabla_{e_i}S))(e_i,e_j)  \quad {\rm using \ (3.4)}
$$
$$
	=\frac 12 \{ (\nabla_{e_i}S)((\nabla_{e_j}J)e_i,Je_j) + (\nabla_{e_i}S)(Je_i,(\nabla_{e_j}J)e_j)
$$
$$
	+S((\nabla_{e_i}(\nabla_{e_j}J))e_i,Je_j)+S(Je_i,(\nabla_{e_i}(\nabla_{e_j}J))e_j)
$$
$$
	- (\nabla_{e_j}S)((\nabla_{e_i}J)e_i,Je_j) - (\nabla_{e_j}S)(Je_i,(\nabla_{e_i}J)e_j)
$$
$$
	- S((\nabla_{e_j}(\nabla_{e_i}J))e_i,Je_j) - S(Je_i,(\nabla_{e_j}(\nabla_{e_i}J))e_j) \} \quad {\rm using \ (2.1)}
$$
$$
	=\frac 12 \{ (\nabla_{e_i}S)((\nabla_{e_j}J)e_i,Je_j) + (\nabla_{e_i}S)(Je_i,(\nabla_{e_j}J)e_j)
$$
$$
	- (\nabla_{e_j}S)((\nabla_{e_i}J)e_i,Je_j) - (\nabla_{e_j}S)(Je_i,(\nabla_{e_i}J)e_j)
	+(\lambda_j-\lambda_i)R(e_i,e_j,e_j,e_i)  \}
$$
and using (3.4) we obtain
$$
	  (\nabla_{e_i}(\nabla_{e_j}S))(e_i,e_j)-
	  (\nabla_{e_j}(\nabla_{e_i}S))(e_i,e_j) = \frac 12 (\lambda_j-\lambda_i)R(e_i,e_j,e_j,e_i) 
$$
$$
	  +\frac 14 (\lambda_j-\lambda_i) \{
	  g((\nabla_{e_i}J)e_j,(\nabla_{e_j}J)e_i) - g((\nabla_{e_i}J)e_i,(\nabla_{e_j}J)e_j) \}.
$$
On the other hand, (2.1) implies
$$
	  (\nabla_{e_i}(\nabla_{e_j}S))(e_i,e_j)-
	  (\nabla_{e_j}(\nabla_{e_i}S))(e_i,e_j) =  (\lambda_j-\lambda_i)R(e_i,e_j,e_j,e_i) 
$$
and hence we find
$$
	\begin{array}{c}
		(\lambda_j-\lambda_i)R(e_i,e_j,e_j,e_i)  \\
		\ds =\frac 12 (\lambda_j-\lambda_i) \{
	  g((\nabla_{e_i}J)e_j,(\nabla_{e_j}J)e_i) - g((\nabla_{e_i}J)e_i,(\nabla_{e_j}J)e_j) \}
	\end{array}      \leqno (3.7)
$$
for all $i,\, j=1,...,2m$. If $e_i \ne e_j,\ Je_j$ we have 
$ R(e_i,Je_j,Je_j,e_i)=R(e_i,e_j,e_j,e_i)$ because of $C=0$. Consequently
from (3.7)  and (1.2) we derive
$$
		(\lambda_j-\lambda_i)R(e_i,e_j,e_j,e_i) 
		 =\frac 12 (\lambda_j-\lambda_i) 
	  g((\nabla_{e_i}J)e_j,(\nabla_{e_j}J)e_i)      \leqno (3.8)
$$
and this is true also for $e_i=e_j$ or $e_i=Je_j$.

From (3.5), (3.6) and (3.8) we obtain
$$
	\sum_{i=1}^{2m} (\lambda_j-\lambda_i) 
	  g((\nabla_{e_i}J)e_j,(\nabla_{e_j}J)e_i) =0     \leqno (3.9)
$$
for any $j=1,...,2m$. Let $ \lambda_1 \le \lambda_2 \le ... \le \lambda_m$. Using (3.9) we find
$$
	\begin{array}{l}
		\quad\ \  \qquad \lambda_i=\lambda_1   \qquad {\rm or} \qquad (\nabla_{e_i}J)e_1=0  \\
		{\rm and}   \qquad   \lambda_i=\lambda_m   \qquad {\rm or} \qquad (\nabla_{e_i}J)e_m=0  
	\end{array}   \leqno (3.10)
$$
for each $i=1,...,m$. If there exists $j$, such that $\lambda_1 < \lambda_j < \lambda_m$, then 
from (3.8) and (3.10) we derive
$$
	R(e_1,e_j,e_j,e_1) = 0, \qquad  R(e_m,e_j,e_j,e_m) = 0 
$$
and because of $C=0$ this implies $\lambda_1=\lambda_m$, which is a contradiction. Consequently 
we have the following two cases:

1) $\lambda_i=\lambda_j$ for all $i,\, j=1,...,m$;

2) $\lambda_i=\lambda$ for $i=1,...,n$, $\lambda_i=\mu$ for $i=n+1,...,m$,
$\lambda \ne \mu$, $1\le n < m$.

In both cases using (3.4) and (3.10), we obtain $\nabla S=0$ in $p$. Consequently
the Ricci tensor is parallel.

If $M$ is irreducible, it is an Einsteinnian manifold and since $M$ is conformal flat, it 
is of constant sectional curvature. Then the Theorem follows from the results in the
end of Preliminaries. 

Let $M$ be a reducible non Einsteinnian manifold. Then we have the case 2) for each 
$p \in M$. Now $M$ is locally a product $M_1 \times M_2$, where $M_1$ and $M_2$ 
are almost K\"ahler manifolds. Let for example $ \dim M_1 \ge 4$. Let $x,\, y$ be
orthogonal unit vectors in a point of $M_1$ and $z$ be a unit vector on $M_2$. 
Because of $C=0$ we have 
$$
	R(x-z,y+Jz, y-Jz, x+z) = 0
$$
or
$$
	R(x,y,y,x) +R(z,Jz,Jz,z) = 0 \ .     \leqno (3.11)
$$
Hence $M_1$ is of constant sectional curvature, say $-c$ and consequently $c \ge 0$.
If $\dim M_2 =2$, it follows from (3.11) that $M_2$ is of constant sectional curvature $c$.
If $\dim M_2 \ge 4$, then $M_2$ is of constant sectional curvature, say $k$ and 
from (3.11) we obtain $k=c$. If $c>0$ this is impossible, because of 
$\dim M_2 \ge 4$ and if $c=0$ $M$ is Einsteinnian, which is a contradiction.

\vspace{0.4in}
\centerline{\large R E F E R E N C E S}

\vspace{0.2in}
\noindent
1. M. B a r r o s. Classes de Chern de las $NK$-variedades. Geometria de las
$AK_2$-varie-

\ \ \ \ \ \ \ \ \ \ 
dades. Tesis doctorales. Universidad de Granada, 1977.

\noindent
2. A. G r a y. Curvature identities for Hermitian and almost Hermitiam manifollds.

\ \ \ \ \ \ \ \ \ \  
{\it T\^ohoku math. J.}, {\bf 28}, 1976, 601-612.

\noindent
3. Z. O l s z a k. A note on almost Kaehler manifolds. {\it Bull. Acad. Polon. Sci.,
S\'er. Sci. 

\ \ \ \ \ \ \ \ \ \ Math. Astr. Phis.}, {\bf 26}, 1978, 139-141.

\noindent
4. K. T a k a m a t s u, \ Y. W a t a n a b e. Classification of a conformally flat
 $K$-space. 
 
\ \ \ \ \ \ \ \ \ \ 
 {\it T\^ohoku math. J.}, {\bf 24}, 1972, 435-440.

\noindent
5. S. T a n n o. 4-dimensional conformally flat K\"ahler manifolds. {\it T\^ohoku math. J.}, {\bf 24},

\ \ \ \ \ \ \ \ \ \ 
   1972, 501-504.

\vspace {0.5cm}
\noindent
{\it Center for mathematics and mechanics \ \ \ \ \ \ \ \ \ \ \ \ \ \ \ \ \ \ \ \ \ \ \ \ \ \ \ \ \ \ \ \ \ \
Received 20.VII.1982

\noindent
1090 Sofia   \ \ \ \ \ \ \ \ \ \ \ \ \ \ \ \ \  P. O. Box 373}

\end{document}